\documentclass[12pt]{amsart}
\usepackage{amsmath, amsthm, url, epsfig, amssymb, setspace, hyperref, mathrsfs}
\usepackage[all]{xypic}

\newtheorem{thm}{Theorem}[section]
\newtheorem{prop}[thm]{Proposition}
\newtheorem{cor}[thm]{Corollary}
\newtheorem{lemma}[thm]{Lemma}

\theoremstyle{definition}
\newtheorem{defe}[thm]{Definition}
\theoremstyle{remark}
\newtheorem{exam}[thm]{Example}

\newtheorem{ntn}[thm]{Convention}

\newtheorem{rem}[thm]{Remark}

\numberwithin{equation}{section}

\def    \blfootnote {\xdef\@thefnmark{}\@footnotetext}

\long\def\symbolfootnote[#1]#2{\begingroup%
\def\thefootnote{\fnsymbol{footnote}}\footnote[#1]{#2}\endgroup}

\def    \H      {\mathrm{H}}
\def    \I      {\mathrm{I}}
\def    \ab     {\mathrm{ab}}

\def    \et     {\mathrm{et}}

\def    \un     {\mathrm{et}}

\def    \id     {\mathrm{id}}

\def    \im     {\mathrm{im}}
\def    \Hom    {\mathrm{Hom}}
\def    \Ext    {\mathrm{Ext}}

\def    \cent   {\mathrm{cent}}
\def    \spec   {\mathrm{spec}}
\def    \red    {\mathrm{red}}
\def    \fppf   {\mathrm{fppf}}

\def    \disc   {\mathrm{disc}}
\def    \Aut    {\mathrm{Aut}}

\def    \coker  {\mathrm{coker}}

\def    \true   {\mathrm{true}}
\def    \Fr     {\mathrm{Fr}}
\def    \SL     {\mathrm{SL}}

\def    \PGL    {\mathrm{PGL}}

\def    \cont  {\mathrm{continuous}}

\def    \st     {\mathrm{st}}
\def    \ad     {\mathrm{ad}}
\def    \ss     {\mathrm{ss}}

\def    \e      {\'{e}}

\def        \Gam            {\Gamma}

\def        \ph             {\varphi}
\def        \GG             {[G,G]}

\def        \uG             {[G,G]^\true}

\def        \uP             {G^{\ab,\st}}
\def        \p              {\Pi_\cent}
\def        \pet            {\Pi_{\et}}
\def        \pdisc          {\Pi_{\disc}}

\def        \uA             {A^\un}

\def        \qz             {\bQ/\bZ}

\def        \d              {\delta}

\def    \Cov    {\mathbf{Cov}}
\def    \bone    {\textbf{1}}

\def    \2Gp    {\textbf{2Gp}}

\def    \bCent  {\mathbf{Cent}}
\def    \bAb    {\mathbf{Ab}}
\def    \bFin   {\mathbf{Fin}}

\def    \2Gr        {\textbf{2Grp}}

\def    \sDMgrst    {\mathsf{DMgrst}}

\def    \cK     {\mathcal{K}}
\def    \cP     {\mathcal{P}}

\def    \cS     {\mathcal{S}}

\def    \cC     {\mathcal{C}}

\def    \cD     {\mathcal{D}}

\def    \cM     {\mathcal{M}}

\def    \bQ     {\mathbb{Q}}

\def    \bC     {\mathbb{C}}
\def    \bZ     {\mathbb{Z}}
\def    \Ga     {\mathbb{G}_a}

\def    \Fq     {\mathbb{F}_q}
\def    \Fp     {\mathbb{F}_p}
\def    \bFq    {\bar{\mathbb{F}}_q}

\def    \bQl    {\overline{\mathbb{Q}}_\ell}
\def    \bA     {\mathbb{A}}

\def    \tG     {\tilde{G}}

\def    \tH     {\tilde{H}}

\def    \tm     {\tilde{m}}

\def    \tT     {\tilde{T}}

\def    \fG     {\mathfrak{G}}

\def    \ra     {\rightarrow}
\def    \inj    {\hookrightarrow}
\def    \onto   {\twoheadrightarrow}

\def    \*      {\times}

\def    \unG    {\underline{G}}

\newcommand         {\Tors}[1]      {#1\!-\!\textbf{Tors}}

\newcommand         {\rar}[1]       {\stackrel{#1}{\longrightarrow}}
\newcommand         {\isom}         {\rar{\simeq}}
\newcommand         {\eqv}          {\rar{\sim}}
\newcommand         {\ilim}         {\mathop{\varprojlim}\limits}

\newcommand         {\commentout}[1]    {}

\begin{document}

\title{Stacky abelianization of algebraic groups}
\author{Masoud Kamgarpour}
\thanks{Author was supported by NSERC PGS D2 - 331456
\newline {\it Address}: {Department of Mathematics, University of Chicago,
Chicago, IL 60637}
 \newline {\it E-mail}: {\url{mkamgarp@math.uchicago.edu}}
 \newline {\it Date}: \today}

\bibliographystyle{alpha}

\begin{abstract} We prove a conjecture of Drinfeld regarding
restriction of central extensions of an algebraic group $G$ to the
commutator subgroup. As an application, we construct the ``true
commutator" of $G$. The quotient of $G$ by the action of the true
commutator is the universal commutative group stack to which $G$ maps.
\end{abstract}
\maketitle

{\small \setcounter{tocdepth}{1} \tableofcontents}


\section{Introduction}
\label{c:introduction} One of the goals of geometric character
theory is constructing sheaves on an algebraic group $G$ over a
finite field $\Fq$ whose ``trace of Frobenius functions" are the
irreducible characters of $G(\Fq)$ ~\cite{Lusztig}, \,
~\cite{Lusztig03}, \newline ~\cite{DB}.
 A {\it one-dimensional character sheaf} is easy to define: it is an
irreducible local system $\cK$ equipped with an isomorphism $m^*\cK
\cong \cK\boxtimes \cK$, where $m:G\times G\ra G$ is the group
multiplication; see for instance, ~\cite{dennis}, \S 4.3.
One-dimensional character sheaves are closely related to central
extensions. Their relationship can be summarized in the following
table
\[
\begin{tabular}{|l|l|}
 \hline
 {\small central ext. $1\ra A\ra \tG\rar{\pi} G\ra 1$}     &   homomorphisms $G(\Fq)\ra A$
 \\
 \hline
 1-dim character sheaves                      &
1-dim characters $G(\Fq)\ra \bQl^\times$\\
 \hline
\end{tabular}
\]
Here $A$ is a finite abelian group. To go from left to right one
takes the trace of Frobenius acting on the stalks of the sheaf (in
the case of central extensions, the sheaf of local sections of
$\pi$). To go from top to bottom, one needs the extra data of a
homomorphism $A\ra \bQl^\times$.

Recall that every one-dimensional character of $G(\Fq)$ is trivial
on the commutator subgroup $[G(\Fq),G(\Fq)]$. There exist, however,
examples of one-dimensional character sheaves whose restriction to
$\GG$ is nontrivial. For instance, we have the central extension
\begin{equation}\label{ex:PGLn}
1\ra \mu_n\ra \SL_n \ra \PGL_n \ra 1.
\end{equation}
Our goal is to construct a commutator for $G$ which is suitable for
doing geometric character theory. To accomplish this, we need to
study the relationship between central extensions of $G$ and those
of $\GG$.

Henceforth, let $G$ denote a connected algebraic group over a
perfect field $k$. We will prove that there exists a pro-\e tale
group scheme $\pet(G)$ such that for every commutative \e tale group
scheme $A$, $\Hom(\pet(G),A)$ equals the set of isomorphism classes
of central extensions of $G$ by $A$.\footnote{This is the geometric
analogue of the following fact: For a perfect group $\Gam$,
$\Hom(\H_2(\Gam,\bZ),A)=\H^2(\Gam,A)$; see for instance,
~\cite{Milnori}.}

\begin{exam} Let $\Ga$ denote the additive group over $\bFq$. Then
\[
\pet(\Ga)=\Hom(\bFq,\bC^\times).
\]
\end{exam}

Restricting central
extensions of $G$ to $\GG$ defines a morphism of group schemes
\begin{equation}
 \pet(i): \pet(\GG)\ra \pet(G).
\end{equation}
The following theorem was conjectured by V. Drinfeld.

\begin{thm} \label{t:finiteness}
The image of $\pet(i)$ is finite.
\end{thm}

Let $\uA$ be the image of $\pet(i)$ and let
\begin{equation}
1\ra \uA\ra \uG\ra \GG\ra 1
\end{equation}
 be the central extension
corresponding to the canonical epimorphism $\pet(\GG)\onto \uA$. We
propose to think of $\GG^\true$ as the ``true commutator" of $G$.
True commutator is characterized by the following properties:
\begin{enumerate}
\item[P1)] $\GG^\true$ is a connected \e tale central extension of
$\GG$.

\item[P2)] The pullback of every \e tale central extension of $G$ to $\GG^\true$ is trivial.

\item[P3)] The commutator map has a lift to $\uG$; that is to say, there exists a
morphism of algebraic varieties $G\times G\ra \uG$ such that the
following diagram commutes
\[
\xymatrix@R=0.3cm@C=0.3cm{
                &         \uG \ar[dd]    \\
  G\times G \ar@{.>}[ur] \ar[dr] _c                \\
                &         \GG,              }
\]
where $c(g,h)=[g,h]:=g^{-1}h^{-1}gh$ for $g,h\in G$.
\end{enumerate}

\begin{rem} Over an algebraically closed field, the fact that P1)-P3)
characterize the true commutator is equivalent to the exactness of
the following complex of profinite groups:
\begin{equation}
\pi_1(G\times G)\ra \pet(\GG)\ra \pet(G).
\end{equation}
\end{rem}

\begin{rem} If $k$  is algebraically closed of characteristic zero and $G$
is reductive, then the true commutator of $G$ is the simply
connected cover of the semisimple group $\GG$, see Example
\ref{ex:reductive}. The fact that the commutator map lifts to
$\GG^\true$ was observed by P. Deligne ~\cite{Del}, \S 2.0.2.
\end{rem}

We now study the quotient of $G$ by the true commutator. Note that
the composition $\GG^\true\onto \GG\inj G$ defines an action of
$\GG^\true$ on $G$ by left translation. The quotient stack
 \begin{equation}
 \uP:=[G/\uG]
 \end{equation}
 is called the {\it stacky abelianization} of $G$. Observe that stabilizers
of the action of $\GG^\true$ on $G$ are isomorphic to $A^\et$. In
particular, stacky abelianization is a Deligne-Mumford stack
~\cite{DM}. We will show, moreover, that it is a strictly
commutative Picard stack ~\cite{DeligneSGA4III}. Roughly speaking,
this means that $\uP$ is equipped with a commutative group
structure. The following theorem states that $\uP$ is the universal
commutative group stack to which $G$ maps.\footnote{For a more
precise formulation, see Theorem \ref{t:UniversalStack}.}

\begin{thm}\label{t:IntrouniversalStack}
Let $\cP$ be a strictly commutative Deligne-Mumford Picard stack,
and let $G\ra \cP$ be a 1-morphism of group stack. Then there exists
a 1-morphism of Picard stacks $\uP\ra \cP$ such that the following
diagram commutes: $ \xymatrix@R=0.3cm@C=0.3cm{
                &         \uP \ar@{.>}[dd]    \\
  G \ar[ur] \ar[dr]                 \\
                &        \cP.              }$
\end{thm}

\subsection{Organization of the text} We start \S \ref{c:centralExt} by proving an algebra-geometric analogue
of the fact that a morphism from a connected space to a finite space
is constant. We then show that every connected algebraic group has a
``universal pro-cover", and use this to define $\pet$. We will see
that over an algebraically closed field, $\pet$ is a quotient of the
Grothendieck fundamental group. We conclude the section by providing
a summary of our calculations of fundamental group schemes. The
details of these calculations appear in Appendix
\ref{a:ComputingFundamentalGroupSchemes}.

In \S \ref{c:trueComandStackyAbe} we study the relationship between
central extensions of $G$ and those of $\GG$. In particular, we
provide examples of central extensions of {\it unipotent} groups
whose restriction to $\GG$ is nontrivial. Moreover, we give a
necessary and sufficient criterion for ``lifting" central extension
of $\GG$ to $G$. The proof of Theorem \ref{t:finiteness} is also
given in this section.

In \S \ref{s:trueComStackAbel}, we study the true commutator and the
stacky abelianization. We will show that the commutator map lifts to
the true commutator and use this to prove that the stacky
abelianization is Picard. We then prove the characterization of the
true commutator and the universal property of the stacky
abelianization (Theorem \ref{t:IntrouniversalStack}).

Throughout the text we assume that the reader is familiar with group
stacks and their relationship to crossed modules. For the
convenience of the reader, we have gathered all the relevant results
in Appendix \ref{a:crossedModandGrStacks}. The only part of this
appendix which appears to be new is the stacky analogue of the First
Isomorphism Theorem (Lemma \ref{l:firstIsoForgrStacks}).

\subsection{Acknowledgments} I wish to express my profound gratitude
to my thesis advisors Vladimir Drinfeld and Mitya Boyarchenko. I
also like to thank  Justin Noel, Larry Breen, Behrang Noohi, and
Travis Schedler for many useful conversations.

\section{Central extensions of algebraic groups}
\label{c:centralExt}

\subsection{Constant morphisms}
Let $\ph:X\ra Y$ be a morphism of schemes over a noetherian
base-scheme $S$.

\begin{defe} \label{d:constantMor}
$\ph$ is {\it constant} if for every scheme $T$ over $S$, the map of
sets $\ph_T: \Hom_S (T,X)\ra \Hom_S (T,Y)$ is constant.
\end{defe}

Let $S'$ be a scheme over $S$.
\begin{lemma}\label{l:baseChangeForConstantMorphims}
\begin{enumerate}
\item[(i)] If $\ph$ is constant, so is $\ph\times_S S'$.
\item[(ii)] Assume $S'\ra S$ is flat and $Y\ra S$ is finite. Then $\ph$ is
constant if and only if $\ph\times_S S'$ is constant.
\end{enumerate}
\end{lemma}
\begin{proof} Statement (i) follows from the fact that for every scheme $T'$
over $S'$,
\[
\Hom_{S'}(T',X\times_S S')=\Hom_S(T',X). \]

Under the conditions of (ii), for every scheme $T$ over $S$, the
base change map
\[
\Hom_S(T,Y)\ra \Hom_{S'}(T\times_S S',Y\times_S S)
\]
 is an
isomorphism (EGA 1, ch. 0, \S 6.2.2). This proves statement (ii).
\end{proof}

\begin{lemma} \label{l:ConstMorEquivDef}
Assume $X$ admits an $S$-section $\sigma: S\ra X$. Then the
following are equivalent:
\begin{enumerate}
\item $\ph$ is constant.
\item $\ph$ has a factorization $X\ra S\ra Y$.
\end{enumerate}
\end{lemma}

\begin{proof}
If $\ph$ is constant then $\ph=(\ph\circ \sigma) \circ h$, where
$h:X\ra S$ is the structure map. Conversely, if $\ph$ has a
factorization $X\ra S\ra Y$ then $\ph_T(\Hom_S(T,X))$ is the set
consisting of one element; namely, the composition $T\ra S\ra Y$.
\end{proof}

\commentout{ Recall that a scheme $A$ over a field $k$ is said to be
finite (resp. \e tale) if the structure morphism $A\ra \spec(k)$ is
finite (resp. finite and \e tale). A finite (resp. \e tale) scheme
over $k$ is isomorphic to $\spec(M)$ where
\begin{equation}\label{eq:prodOfAlgebras}
M\cong \prod_{i}^n M_i,
\end{equation}
 and each $M_i$ is a finite
dimensional $k$-algebra (resp. an \e tale $k$-algebra).
 We say that $A$ is {\it discrete}\footnote{This terminology is not
standard.} if $M_i\cong k$ for every $i$. The following proposition
will be used repeatedly. }

\begin{prop} \label{p:ConnectedToFiniteIsConstant}
Let $k$ be a perfect field, $A$ a finite scheme over $k$, and $X$ a
connected reduced scheme over $k$ which admits a $k$-section. Then
every morphism $\ph:X\ra A$ of $k$-schemes is constant.
\end{prop}

\begin{rem} \label{r:connectedGeometric}
As $k$ is perfect, $X$ is reduced if and only if it is geometrically
reduced (EGA IV, part 2, prop. 4.6.1). As $X$ admits a $k$-section,
it is connected if and only if it is geometrically connected (loc.
cit., cor. 4.5.14).
\end{rem}

\begin{proof}[Proof of Proposition \ref{p:ConnectedToFiniteIsConstant}]
In view of Lemma \ref{l:baseChangeForConstantMorphims} and Remark
\ref{r:connectedGeometric}, we may assume that $k$ is algebraically
closed. We may suppose, moreover, that $X$ is affine.

Let $X=\spec(R)$ and $A=\spec(M)$, where $R$ and $M$ are
$k$-algebras. Let $\Phi: M\ra R$ be the $k$-algebra homomorphism
corresponding to $\ph$. Let $N\subseteq R$ be the image of $\Phi$.
Since $R$ has no nontrivial idempotents, $N$ is a finite dimensional
local $k$-algebra. As $R$ is reduced, $N$ is isomorphic to $k$. By
Lemma \ref{l:ConstMorEquivDef}, $\ph$ is constant.
\end{proof}

\commentout{
\begin{cor} Let $G$ be a connected algebraic group and let $A$ be an \e tale normal sub-group-scheme
of $G$. Then $A$ is central.
\end{cor}
\begin{proof} Let $G'$ be a connected component of
$G\times A$. Then $G'$ is (non-canonically) isomorphic to $G$; in
particular, it is geometrically reduced and connected. By
Proposition \ref{p:ConnectedToFiniteIsConstant}, the morphism
$G'\times A\ra A$ defined by $g'.a:=g'^{-1}ag$ is trivial. Hence the
corresponding morphism $G\times A \ra A$ is constant.
\end{proof}

\begin{rem} Let $G$ be a connected algebraic group and let $A$ be an \e tale normal subgroup scheme
of $G$. Using the above proposition one can easily show that $A$ is
central.
\end{rem}
}

\begin{rem} (M. Boyarchenko) Let $X$ be as in Proposition \ref{p:ConnectedToFiniteIsConstant},
$S$ be a {\it reduced} scheme over $k$, and $A$ be a finite scheme
over $S$. Then one can show that every $S$-morphism
$X\times_{\spec(k)} S\ra A$ is constant. On the other hand, if $S$
is not reduced, then this results fails. For example, let $T$ be a
$k$-algebra containing a nonzero element $u$ satisfying $u^d=0$, for
some positive integer $d$. Let $X=\bA_k^1$, $S=\spec(T)$, and
$A=\spec(T[\epsilon ]/\epsilon^d)$. The map $\epsilon \mapsto u.x$
extends to a $T$-algebra homomorphism $T[\epsilon ]/\epsilon^d \ra
T[x]$. The corresponding morphism $X\times_{\spec(k)} S\ra A$ of
schemes over $S$ is not constant.
\end{rem}

\subsection{Fundamental group schemes of algebraic groups}
 \label{s:FundamentalGrouPSchemes} Let $k$ be a perfect field.
 Let $\fG$ denote the category of
connected algebraic groups over $k$.

\begin{lemma} \label{l:fibreProductExists}
Fiber products exists in $\fG$.
\end{lemma}

\begin{proof} Given a diagram $G_1\ra G\leftarrow G_2$ in $\fG$, its
fiber product equals $(G_1\times_G G_2)_\red^0$: the reduced neutral
connected component of the scheme-theoretic fibre product
$G_1\times_G G_2$.
\end{proof}

\begin{defe} A morphism $\tG\ra G$ in $\fG$
 is an {\it isogeny} if it is surjective and its kernel is
finite. (Of course, the kernel is computed in the category of {\it
all} group schemes.)
\end{defe}

\begin{defe}\label{d:centralCover}
A {\it group cover} of $G$ is a connected algebraic group $\tG$
equipped with an isogeny $\tG\ra G$. The group $\tG$ is a {\it
central cover} if the kernel of $\tG\ra G$ is central.\footnote{It
is easy to show that an {\it \e tale} group cover is central.}
\end{defe}

Let $\Cov(G)$ denote the category of coverings of $G$. (An arrow in
$\Cov(G)$ is a morphism $\tG\ra \tG'$ of group schemes over $G$.) It
is clear that $\Cov(G)$ is essentially small.

\begin{lemma}\label{l:coveringsIsadirectedSet}
$\Cov(G)$ is anti-equivalent to a partially ordered directed set
$I=I(G)$.
\end{lemma}

\begin{proof} Proposition \ref{p:ConnectedToFiniteIsConstant} shows
that there is at most one morphism between two coverings of $G$.
Thus, $\Cov(G)$ is a partially ordered set. The supremum of two
element of $\Cov(G)$ is given by their fibre product (Lemma
\ref{l:fibreProductExists}).
\end{proof}

For each $i\in I(G)$, let $G_i\ra G$ be the corresponding object of
$\Cov(G)$, and let $A_i:=\ker(G_i\ra G).$ The kernels $A_i$'s form
an inverse system: there is a morphism $A_i\ra A_j$ if and only if
there is a morphism $G_i\ra G_j$ in $\Cov(G)$. The set $\I(G)$ has
the following subsets:
\begin{equation}
\I_{\cent}(G)\supseteq \I_{\et}(G)\supseteq \I_\disc(G),
\end{equation}
corresponding to group covers of $G$ with central, \e tale, and
discrete kernel. (A discrete group scheme is a finite \e tale scheme
on which the Galois group acts trivially.)

Every element of $\I$ dominated by an element of $\I_\cent$ belongs
to $\I_\cent$. Furthermore, the supremum of two elements of
$\I_\cent$ belongs to $\I_\cent$. These facts remain true if
$\I_\cent$ is replaced by $\I_\et$ or $\I_\disc$.

\begin{ntn} Set
\[\label{q:versionsoffundgr} \p(G):=\ilim_{i\in I_\cent(G)}A_i,
\]
\[
\pet(G):=\ilim_{i\in I_\et(G)}A_i,
\]
\[\pdisc(G):=\ilim_{i\in
I_\disc(G)}A_i.
\]
We refer to these profinite group schemes as the {\it fundamental
group schemes} of $G$.
\end{ntn}

\begin{rem} \label{r:petisquotientOfp}
The group scheme $\pet(G)$ is the maximal pro\e tale quotient of
$\p(G)$. The group $\pdisc(G)$ is the maximal quotient of $\pet(G)$
on which the absolute Galois group of $k$ acts trivially.
\end{rem}

The following table summarizes our computations of fundamental group
schemes of certain connected algebraic groups. See Appendix
\ref{a:ComputingFundamentalGroupSchemes} for the details.
\begin{ntn} For a semisimple group $G$, $\pi_1^\ss(G)$ denotes the
weight lattice modulo the root lattice ~\cite{Borel}, \S 24.1. For
an algebraic group $G$, $\pi_1(G)$ denotes the algebraic fundamental
group of $G$ ~\cite{GR}. In characteristic zero,
$\pi_1^\ss(G)=\pi_1(G)$.
\end{ntn}

\[
\begin{tabular}{|l|l|l|}
 \hline
{\small $G$  defined over $\bC$}                 &
$\pet(G)=\pi_1(G)$
 \\
 \hline
 $G$ semisimple                   & $\pet(G)=\pi_1^\ss (G)$
 \\
 \hline
 $G$ commutative defined over $\Fq$   &   $\pdisc(G)=G(\Fq)$\\
 \hline
 $G=\Ga$ additive group over $\overline{\Fq}$               &   $\pet(G)=\Hom(k,\qz)$\\
 \hline
\end{tabular}
\]

\subsection{Classifying central extensions}
Let $A$ be a finite commutative group scheme over $k$. Let
\[
1\ra A\ra \tG\ra G\ra 1
\]
be a central extension. Then $(\tG)_\red^0\ra G$ is a central cover.
Let $f_{\tG}:\p(G)\ra A$ denote the composition
\[
\p(G)\onto \ker((\tG)_\red^0\ra G)\inj A.
\]
\begin{prop} \label{p:petClassifiesCentral}
The map $\tG\mapsto f_{\tG}$ defines a canonical isomorphism
\begin{equation}
\label{eq:Hom=H^2}
 \H^2(G,A)\isom \Hom(\p(G),A).
\end{equation}
If $A$ is \e tale \emph{(}resp. discrete\emph{)}, then we can
replace $\p(G)$ by $\pet(G)$ \emph{(}resp. $\pdisc(G)$\emph{)}.
\end{prop}

\begin{proof} Let $\ph \in \Hom(\pet(G),A)$. Let $B:=\im(A)$. The
group $B$ is a finite quotient of $\p(G)$; thus, $B\cong A_i$ for
some $i\in I_\cent$. Let $G^\ph$ be the pushforward of the central
extension
\[
1\ra A_i\ra G_i\ra G\ra 1
\]
along the morphism $A_i\isom B\inj A$. One checks that $\ph\mapsto
G^\ph$ is the inverse of $\tG\ra f_{\tG}$.
\end{proof}

\subsection{$\pet$ is a quotient of $\pi_1$} In this section, we assume $k$ is
algebraically closed. Let $\pi_1(G):=\pi_1(G,e)$ denote the
fundamental group of $G$ in the sense of ~\cite{GR}. For every
abelian group $A$ one has a homomorphism
\begin{equation}
\xymatrix{ \Hom(\pet(G),A) \ar@{=}[r] \ar[d]^{f_A} &
\{ \textrm{central extensions of $G$ by $A$} \} \\
\Hom(\pi_1(G),A)  \ar@{=}[r]                       &
\{\textrm{$A$-torsors on $G$ trivialized over $e$} \}. }
\end{equation}
Since $f_A$ is functorial in $A$, it comes from a homomorphism
$f:\pi_1(G)\ra \pet(G)$.

\begin{lemma}\label{l:petQuitient}
\begin{enumerate}
\item[(i)] $f_A$ is injective.
\item[(ii)] $\ph: \pi_1(G)\ra A$ belongs to the image of $f_A$ if and only if the
diagram
\begin{equation}
\xymatrix{
 \pi_1(G\times G)\ar[r]^{m_*}   \ar[d]  &     \pi_1(G) \ar[r]        & A \\
 \pi_1(G)\times \pi_1(G) \ar[d]         &                       &  \\
 A\times A \ar[uurr]                     &                       &
 }
\end{equation}
commutes.
\end{enumerate}
\end{lemma}
\begin{proof}
A based $A$-torsor $\tG\ra G$ is a central extension if and only if
$m^*\tG$ is isomorphic to  $\tG\boxtimes \tG$ (as $A$-torsors on
$G\times G$). This proves (i).

To prove (ii), let $\tG\ra G$ denote the based $A$-torsor
corresponding to $\ph$. One checks that the commutativity of the
diagram is equivalent to the existence of an ($A\times
A$)-equivariant morphism of based schemes $\tm:\tG\times \tG\ra \tG$
such that the diagram
\[
\xymatrix{ \tG\times \tG \ar[r]\ar[d]   & \tG \ar[d] \\
            G\times G   \ar[r]          &   G
}
\]
commutes. Proposition \ref{p:ConnectedToFiniteIsConstant} then
implies that $\tm$ satisfies the axioms defining a group.
\end{proof}

 The following is a reformulation Lemma \ref{l:petQuitient}.

\begin{cor} \label{c:petQuotpi1}
 $\pet(G)$ equals the coequalizer of the following
homomorphisms:
 \begin{gather} \label{dia:coequalizer}
 \pi_1(G\times G)^\ab \rar{m_*}
\pi_1(G)^\ab \\
\notag \xymatrix{\pi_1(G\times G)^\ab \ar[r] \ar@/ _1.7pc/[rr] &
\notag \pi_1 (G)^\ab \times \pi_1 (G)^\ab \ar[r]^{\quad \quad+} &
\pi_1 (G)^\ab.}\\ \notag
\end{gather}
\end{cor}

\begin{rem}
Corollary \ref{c:petQuotpi1} remains valid if one replaces $\pet(G)$
with $\p(G)$ and $\pi_1(G)$ with the fundamental group scheme of $G$
defined ~\cite{Nori}.
\end{rem}

\section{Relationship between central covers of $G$ and $\GG$}
\label{c:trueComandStackyAbe}
\subsection{Restricting central covers of $G$ to $\GG$}
\label{s:pitfallOfCommutator} Let $k$ be a perfect field, $G$ a
connected algebraic group over $k$, and $\tG\ra G$ a central cover
(Definition \ref{d:centralCover}). It is clear that $(\tG \times_G
\GG)_\red^0$ is a central cover of $\GG$.

\begin{ntn} \label{n:RestrictionOfCovers}
We call $(\tG \times_G \GG)_\red^0$ the {\it restriction} of $\tG\ra
G$ to $\GG$. Restricting central covers defines a homomorphism
\[
\p(\GG)\ra \p(G).
\]
\end{ntn}

In the introduction, we mentioned that if $\tG\ra G$ is a nontrivial
central cover of a {\it semisimple} non-simply-connected group $G$,
then the restriction of $\tG$ to $\GG$ is nontrivial. We now give
examples of covers of {\it unipotent} groups whose restrictions to
$\GG$ is nontrivial.

\begin{exam}Let $\tG$ be a connected noncommutative
unipotent algebraic group over a field of positive characteristic.
Let $A$ be an arbitrary finite subgroup of the center $Z(\tG)$ such
that $A\cap [\tG,\tG]\neq \{1\}$.\footnote{For example, $A$ can be
any nontrivial finite subgroup of $C(\tG):=Z(\tG)\cap [\tG,\tG]$.
Note that $\dim(C(\tG))>0$ since the last nonzero term of the lower
central series of $\tG$ is a subgroup of $C(\tG)$. In particular,
$C(\tG)$ has many finite subgroups.} Let $G:=\tG/A$. Then
$\tG\rar{\pi} G$ is a central cover of $G$ whose restriction to
$\GG$ is nontrivial.
\end{exam}

\begin{ntn} For a central cover $\tG\ra G$, let $d(\tG)$ denote the
degree of the central cover $(\tG \times_G \GG)_\red^0\ra \GG$.
\end{ntn}

\begin{thm} \label{t:finiteTheorem}
There exists a constant $C$, depending only on $G$, such that
$d(\tG)<C$ for all central covers $\tG\ra G$.
\end{thm}

In view of Remark \ref{r:petisquotientOfp}, the corresponding
results for $\pet$ and $\pdisc$ follow immediately. In particular,
the image of the map
\[
\pet(\GG)\ra \pet(G)
\]
is finite (Theorem \ref{t:finiteness}).

To prove Theorem \ref{t:finiteTheorem}, we need a lemma. For every
positive integer $n$, let $c^n:G^{2n}\ra \GG$ denote the map
\begin{equation}\label{e:powerOfCommutator}
(g_1,g_2,g_3,g_4,..., g_{2n-1},g_{2n})\mapsto
[g_1,g_2][g_3,g_4]...[g_{2n-1},g_{2n}].
\end{equation}

\begin{lemma} \label{l:powerofcommutatorsurjective}
There exists a positive integer $n$ such that $c^n$ is
surjective.\footnote{In fact, one can take $n=2\dim(G)$.}
\end{lemma}
\begin{proof} See, for instance, ~\cite{MI},
cor. 11.13.
\end{proof}

 The following remark is a key observation used in the proof of Theorem
 \ref{t:finiteTheorem}.

\begin{rem} \label{r:commutatorMapLifts} Let $\tG$ be a central
cover of $G$ and let $H$ be its restriction to $\GG$. The commutator
map $\tG\times \tG\ra \tG$ descends to a morphism of algebraic
varieties $G\times G\ra H$. In this case, we say that the commutator
map of $G$ lifts to $H$.
\end{rem}

\begin{proof}[Proof of Theorem \ref{t:finiteTheorem}]
Choose $n$ large enough so that $c^n$ is surjective. Remark
\ref{r:commutatorMapLifts} implies that $c^n$ has a lift to $H$;
that is to say, there exists a morphism of algebraic varieties
$G^{2n}\ra H$ such that the following diagram commutes:
\begin{equation}
\xymatrix{
                        &       H \ar[dd]^{\textrm{degree $=d(\tG)$}}\\
 G^{2n}  \ar@{.>}[ur] \ar[dr]^{c^n}    &  \\
                        &       \GG\\
}
\end{equation}
As $c^n$ is surjective, there exists a closed subvariety $X\subseteq
G^{2n}$ such that the generic fiber of the morphism $c^n|_X: X\ra
[G,G]$ is finite. The number $d(\tG)$ divides the degree of
$c^n|_X$; hence, it is bounded.
\end{proof}

\subsection{Lifting central covers from $\GG$ to $G$}
\label{ss:liftingCentralExt}

 We keep the
notation and conventions of the previous section. Let $H\ra \GG$ be
a central cover. The following result is a converse of Remark
\ref{r:commutatorMapLifts}.

\begin{prop}\label{p:LiftingCentExt} Suppose the commutator map of $G$ lifts
to $H$. Then there exists a central cover $\tG\ra G$ whose
restriction to $\GG$ is isomorphic to $H$.
\end{prop}

The proof of this proposition will take up the rest of this section.
Let $H$ be a central cover of $\GG$ which has a lift of the
commutator map. Let $\delta: H\ra G$ denote the composition $H\onto
\GG\inj G$.

\begin{lemma}\label{l:LiftoFcommImpliesStrictlyStable}
One can endow  $H\rar{\delta} G$ with a structure of
 is a strictly stable crossed module (Definition
\ref{d:crossMod}).
\end{lemma}

\begin{proof}Let $\{-,-\}:G\times G\ra H$ denote the lift of the commutator.
We may assume that $\{1,1\}=1$. Define a morphism of varieties
$G\times H\ra H$ by
\[
(g,h)\mapsto h^g:=h \{\delta(h),g\}.
\]
 Using Proposition
\ref{p:ConnectedToFiniteIsConstant}, it is easy to show that this
morphism defines an action of $G$ on $H$ making $H\rar{\delta} G$
into a strictly stable crossed module.\footnote{Note that we do not
use the fact that $H$ is connected.}
\end{proof}

In the proof of Proposition \ref{p:LiftingCentExt}, we employ a
theorem about extension of sheaves. To state this result, we need
some notation. Let $k_\fppf$ denote category of schemes over $k$
equipped with the topology generated by faithfully flat morphisms of
finite type. In what follows, $\Ext$ will denote the extension of
abelian sheaves in $k_\fppf$. Let $\bAb$ (resp. $\bFin_k$) denote
the category of finite abelian group (resp. the category of finite
commutative group schemes over $k$). Let $J$ be a connected
commutative algebraic group over $k$.

\begin{thm} \label{t:effaceable}
$\Ext^2(J,-):\bFin_k\ra \bAb$ is an effaceable functor; that is to
say, given $A\in \bFin_k$ and $\alpha\in \Ext(J,A)$, there exists a
monomorphism $A\inj B$ in $\bFin_k$, such that the image of $\alpha$
under the map $\Ext^2(J,A)\ra \Ext^2(J,B)$ is zero.
\end{thm}

\begin{rem}
This theorem was formulated and proved by V. Drinfeld (unpublished).
It is closely related to the fact that $\Ext^2(G,\qz)$ vanishes. For
a proof of this vanishing result see ~\cite{MityaThesis}, lem.
3.2.2.
\end{rem}

\begin{ntn} All Picard stacks considered are assumed to be strictly
commutative, see Definition \ref{d:grCat}.
\end{ntn}

\begin{proof}[Proof of Proposition \ref{p:LiftingCentExt}]

\underline{Step 1}: Let $H$ be a central cover of $\GG$ equipped
with a lift of the commutator map. By Lemma
\ref{l:LiftoFcommImpliesStrictlyStable}, $\d: H\ra G$ is a strictly
stable crossed module. As explained in \S \S
\ref{s:relationShipCrosModGrCat}-\ref{s:grStPicardSt}, this implies
that the quotient stack $\cP$ of $G$ by the action of $H$ is a
Deligne-Mumford
 Picard stack. Set $A:=\ker(\d)$. Note that
 \begin{equation}
 \pi_0(\cP)=\coker(\d)=G^\ab, \quad \quad \textrm{and} \quad \quad
 \pi_1(\cP)=\ker(\d)=A.
 \end{equation}

\underline{Step 2}: It follows from ~\cite{DeligneSGA4III}, prop.
1.4.15, that the set of isomorphism classes of Picard stacks with
$\pi_0=G^\ab$ and $\pi_1=A$ equals $\Ext^2(G^\ab,A)$. Let $\alpha\in
\Ext^2(G^\ab,A)$ denote a representative for the isomorphism class
of $\cP$. By Theorem \ref{t:effaceable}, there exists a monomorphism
$\kappa: A\inj B$ of finite commutative group schemes, such that the
image of $\alpha$ under the induced morphism $\Ext^2(G^\ab,A)\ra
\Ext^2(G^\ab,B)$ is zero. Let $H'$ be the unique central extension
of $\GG$ for which the following diagram is commutative
\footnote{$H'$ is the pushforward of $H$ with respect to $\kappa$.}
 \begin{equation}
\xymatrix{
    1\ar[r]               & A\ar[r] \ar[d]^\kappa & H \ar[r] \ar[d] &  \GG\ar[r] \ar[d] & 1 \\
    1\ar[r]               & B\ar[r]        & H' \ar[r]        &  \GG \ar[r]         &
    1.
}
\end{equation}

\underline{Step 3}: The composition $G\times G\rar{\{-,-\}} H\ra H'$
endows $H'$ with a lift of the commutator map. Applying Lemma
\ref{l:LiftoFcommImpliesStrictlyStable}, we conclude that $H'\ra G$
is a strictly stable crossed module. Let $\cP'$ denote the
corresponding quotient stack. The class of $\cP'$ in
$\Ext^2(G^\ab,B)$ is trivial; therefore,
\[
\cP'\cong G/H\times \Tors{B}.
\]
Here $G/H$ is is the discrete Picard stack defined by the
commutative algebraic group $G/H$ (Example \ref{ex:unG}), and
$\Tors{B}$ denotes the (Deligne-Mumford Picard) stack of $B$-torsors
on $k_\fppf$. The composition
\[
 G \ra \cP':=[G/H']\eqv G \times \Tors{B} \ra
 \Tors{B},
\]
defines a 1-morphism of gr-stacks $G\ra \Tors{B}$.

\underline{Step 4}: Let $\Hom(G,\Tors{B})$ (resp. $\bCent(G,B)$)
denote the Picard stack of 1-morphisms of gr-stacks $G\ra \Tors{B}$
(resp. central extensions of $G$ by $B$). According to ~\cite{SGA7},
\S 1.1, these two Picard stacks are naturally equivalent. Therefore,
the 1-morphism of gr-stacks $G\ra \Tors{B}$ defined in the previous
step, gives rise to a central extension
\[ 1\ra
B\ra G'\ra G\ra 1.
\]
This central extension is a lift of the central extension
corresponding to the 1-morphism of gr-stacks
\[
\GG \ra \Tors{A} \rar{\kappa} \Tors{B}.
\]
 It follows that the
restriction of $G'$ to $\GG$ equals $H'$. Moreover,
$(H')_\red^0\cong H$, as required.
\end{proof}

\section{True commutator and stacky abelianization}
\label{s:trueComStackAbel}

\subsection{Stacky abelianization is Picard}
\label{ss:commutatorMapLiftsTrueCom}
 Let $k$ be an algebraically
closed field and $G$ be a connected algebraic group over $k$. Let
$A^\et$ be the image of $\pet(\GG)\ra \pet(G)$. By Theorem
\ref{t:finiteTheorem}, $A^\et$ is finite. By Proposition
\ref{p:petClassifiesCentral}, we obtain a central extension
\[
1\ra A^\et\ra \GG^\true\ra \GG\ra 1.
\]
We call $\GG^\true$ the {\it true commutator} of $G$. We refer to
the quotient stack of $G$ by the action of $\GG^\true$ as the {\it
stacky abelianization} of $G$ and denote it by $\uP$.

\begin{exam}
\label{ex:reductive} Let $G$ be a connected {\it reductive} group
over an algebraically closed field of characteristic zero. Let
$G^\ad:=G/Z(G)$ be the associated adjoint semisimple group. The
group cover $\GG \ra G^\ad$ defines an injection $\pet(\GG)\inj
\pet(G^\ad)$. It follows that the natural morphism $\pet(\GG)\ra
\pet(G)$ is also injective. As we will see in \S \ref{s:petSSGroup},
$\pet(\GG)=\pi_1^\ss(\GG)$. Therefore $\GG^\true$ identifies with
the simply connected cover of $\GG$.
\end{exam}

\begin{lemma}
\label{l:truecommutatorcomesFromG} The true commutator is the
restriction of a central cover of $G$ (Conventions
\ref{n:RestrictionOfCovers}).
\end{lemma}

To prove this lemma, we need an easy result from the theory of
profinite groups whose proof we omit.

\begin{lemma}\label{l:profiniteGroups}
Let $A\inj C$ be an inclusion of a finite group into a profinite
abelian group. Then there exists an epimorphism $C\onto B$, where
$B$ is a finite group, such that the composition $A \inj C \onto B$
is an injection.
\end{lemma}

\begin{proof}[Proof of Lemma \ref{l:truecommutatorcomesFromG}]
By Lemma \ref{l:profiniteGroups}, there exists a finite quotient $B$
of $\pet(G)$ such that the composition $A^\et\inj \pet(G)\onto B$ is
an injective morphism. Let $\tG$ be the central cover of $G$
corresponding to $\pet(G)\onto B$. The restriction of $\tG$ to $\GG$
is isomorphic to $\GG^\true$.
\end{proof}

\begin{cor} \label{c:commutatorMapLiftstoTrue}
The commutator map lifts to $\GG^\true$.
\end{cor}

\begin{proof} This follows from Lemma \ref{l:truecommutatorcomesFromG}
and Remark \ref{r:commutatorMapLifts}.
\end{proof}

\begin{cor} $\uP$ is a (strictly commutative) Picard stack.
\end{cor}
\begin{proof} By Corollary \ref{c:commutatorMapLiftstoTrue} and
Lemma \ref{l:LiftoFcommImpliesStrictlyStable}, $\uG\ra G$ is a
strictly stable crossed module. Therefore, the corresponding
quotient is a Picard stack.
\end{proof}

\subsection{Characterization of the true commutator}
We keep the notation and conventions of the previous section. Let
$H$ be an \e tale central cover of $\GG$.

\begin{lemma}\label{l:equivalenceOfTrivialRestriction}
The following are equivalent:
\begin{enumerate}
\item[(i)] The pullback
of every central extension of $G$ - by an \e tale group scheme - to
$H$ is trivial.
\item[(ii)] For every \e tale central cover $\tG\rar{\pi} G$,
we have a morphism $H \ra \pi^{-1}(\GG)_\red^0$ of algebraic groups
over $\GG$.
\end{enumerate}
\end{lemma}

\begin{cor} \label{c:trueComMapsToLifts}
Assume that the commutator map lifts to $H$. Then we have a morphism
$\GG^\true\ra H$ of algebraic groups over $\GG$.
\end{cor}
\begin{proof} By Proposition \ref{p:LiftingCentExt}, $H$ is a
restriction of an \e tale central cover of $G$. By the previous
lemma, we obtain a morphism $\GG^\true \ra H$.
\end{proof}

\begin{prop}\label{p:characterizingTruecommutator}
The true commutator is the unique, up to isomorphism, \e tale
central cover of $\GG$ satisfying the following properties:
\begin{enumerate}
\item [P1)] The pullback of every central extension of $G$ - by an
\e tale group scheme - to $\GG^\true$ is trivial.
\item [P2)] The commutator map lifts to $\GG^\true$.
\end{enumerate}
\end{prop}

\begin{proof}[Proof of Proposition
\ref{p:characterizingTruecommutator}] Let $H$ be an \e tale central
cover of $\GG$. Suppose the true commutator lifts to $H$. Then by
Corollary \ref{c:trueComMapsToLifts}, we have a morphism $\GG^\true
\ra H$. Assume, furthermore, that the pullback of every central
extension of $G$ - by an \e tale group scheme - to $H$ is trivial.
In view of the fact that $\GG^\true$ is a restriction of an \e tale
cover of $G$, Lemma \ref{l:equivalenceOfTrivialRestriction} provides
us with a morphism $H\ra \GG^\true$.
\end{proof}

\subsection{Universal property of stacky abelianization}
\begin{ntn}
All stacks discussed below are Deligne-Mumford stacks. The
2-category of Deligne-Mumford gr-stacks is denoted by $\sDMgrst$.
All Picard stacks are assumed to be strictly commutative.
\end{ntn}

Let $\Phi:G\ra \uP$ denote the canonical 1-morphism of gr-stacks.
\begin{thm}\label{t:UniversalStack}
Let $\cP$ be a Picard stack over $k$. Then composition with $\Phi$
defines an equivalence of Picard groupoids
\[
\Hom_{\sDMgrst}(\uP,\cP)\eqv \Hom_{\sDMgrst}(G,\cP).
\]

\end{thm}

\begin{lemma} \label{l:AutoOfGrFunctors}
Let $\cC$ and $\cC'$ be gr-stacks over $k$. Assume that $\pi_0(\cC)$
is representable by a connected algebraic group over $k$. Then, the
groupoid $\Hom_{\sDMgrst}(\cC,\cC')$ is discrete (i.e., the objects
have no nontrivial automorphisms).
\end{lemma}
\begin{proof}
Let $F$ be a 1-morphism of stacks $\cC\ra \cC'$. By Remark
\ref{r:morphismOfGrpi0pi1}, a 2-morphism $F\implies F$ defines a
morphism of schemes $\epsilon: \pi_0(\cC)\ra \pi_1(\cC')$. As $\cC'$
is a Deligne-Mumford stack, $\pi_1(\cC')$ is finite. By Proposition
\ref{p:ConnectedToFiniteIsConstant}, $\epsilon$ is constant.
\end{proof}

\begin{proof}[Proof of Theorem \ref{t:UniversalStack}] In view of the above lemma, it
is enough to show that every 1-morphism of gr-stacks $F:\unG\ra \cP$
has a canonical factorization
\begin{equation}\label{eq:factorizationForUniversal}
G \ra \uP\rar{\Phi} \cP.
\end{equation}
 By the First Isomorphism Theorem for gr-stacks (Lemma \ref{l:FirstIsoForPicardStacks}) and
the remarks following it, there exists a strictly stable crossed
module of algebraic groups $\delta: \tH\ra G$ such that $F$ has a
factorization $\unG\onto [G/\tH]\inj \cP$. Let $H:=\delta(\tH)$.
Then $H$ is a subgroup of $\GG$ containing $\GG$, and $\tH$ is a
central extension of $H$ by an \e tale group scheme. Without loss of
generality, we may assume $H=\GG$.

By Remark \ref{c:trueComMapsToLifts}, we have a morphism of
algebraic groups $\GG^\true\ra H$, which in turn, defines a morphism
of strictly stable crossed modules $(\uG\ra G)\ra (H\ra G)$. Hence,
we obtain a 1-morphism of Picard stacks $\uP=[G/\uG]\ra [G/\tH]\isom
\cP$, providing the required factorization.
\end{proof}

\appendix
\section{Crossed modules and gr-stacks}
\label{a:crossedModandGrStacks}

\subsection{Gr-categories}
\label{s:grCategories}

\begin{ntn} All categories we consider are essentially small.
A monoidal category is  denoted by $(\cM,\otimes, \bone)$. In other
words, we suppress the associativity and unit constraints
~\cite{MacLane}. Occasionally, we suppress $\otimes$ and $\bone$ as
well. A monoidal category is {\it strict} if the associativity and
unit constraints are trivial. With the usual abuse of notation,
$x\in \cM$ means $x$ is an object of $\cM$. We  denote the set of
isomorphism classes of objects of $\cM$ with $\pi_0(\cM)$. For every
$x\in \cM$, $\pi_1(\cM,x)$ denotes the abelian group $\Aut_\cM(x)$.
We set$\pi_1(\cM):=\pi_1(\cM,\bone)$.
\end{ntn}

\label{ss:grcategories}
\begin{defe} \label{d:grCat}
A {\it gr-category} is a monoidal groupoid all of whose objects have
a weak inverse; that is to say, for every $x\in \cM$ there exists
$y\in \cM$ such that $x\otimes y\cong y\otimes x\cong \bone$
~\cite{Sinh}, \S1. $\cM$ is a {\it strict} gr-category if it is a
strict monoidal groupoid such that for every $x\in \cM$ there exists
$y\in \cM$ satisfying $x\otimes y=y\otimes x=\bone$. A (strictly
commutative) {\it Picard groupoid} is a (strictly) symmetric
gr-category ~\cite{DeligneSGA4III}, \S 1.4.
\end{defe}

\begin{exam} \label{ex:unG}
Let $G$ be a group. The discrete groupoid whose set of objects
equals $G$ is a gr-category which, by an abuse of notation, is also
denoted by $G$. If $G$ is commutative, the corresponding gr-category
is a strictly commutative Picard groupoid.
\end{exam}

\begin{rem} \label{r:AutoInGrCat}
Let $\cM$ be a gr-category. For every $x\in \cM$, the map $u\mapsto
u\otimes \id_x$ defines an isomorphism $\pi_1(\cM, \bone) \isom
\pi_1(\cM,x)$.
\end{rem}

\begin{defe} A {\it 1-morphism} of gr-categories (resp. Picard
groupoids) is a monoidal functor (resp. symmetric monoidal functor)
$F: \cM\ra \cM'$. $F$ is a monomorphism (resp. epimorphism, resp.
isomorphism)  if it is fully faithful (resp. essentially surjective,
resp. an equivalence). Gr-categories and Picard groupoids form a
2-category where 2-morphism are monoidal natural transformations.
\end{defe}

\begin{rem}\label{r:morphismOfGrpi0pi1}
Let $\cM$ and $\cM'$ be gr-categories. Let $F$ be a (not necessarily
monoidal) functor $\cM\ra \cM'$, and let $\eta$ be a natural
transformation $F \implies F$. For every $m\in \cM$, we obtain an
element of $\Aut_{\cM'}(F(m))$ which, by the identification of
Remark \ref{r:AutoInGrCat}, gives us an element of $\pi_1(\cM')$.
One checks that this element depends only on the isomorphism class
of $m$ in $\cM$. Therefore, we obtain a well-defined map $\eta_0:
\pi_0(\cM)\ra \pi_1(\cM')$. Note that if $\eta_0$ is trivial, then
so is $\eta$.
\end{rem}

\subsection{Crossed modules}
\begin{defe} \label{d:crossMod}
A (right) {\it crossed module} $\fG$ is the data consisting of a
group homomorphism $\delta: H\ra G$ and a right action of $G$ on
$H$, denoted by $h\mapsto h^g$, such that for every $g\in G$ and
$h_1,h_2\in H$,
\[
(h_2)^{\delta(h_1)}=h_1^{-1}h_2h_1, \quad \quad
\delta(h^g)=g^{-1}\delta(h)g.
\]
$\fG$ is {\it strictly stable} if there exists a map $\{-,-\}:
G\times G\ra H$ such that for every $g,g_0,g_1,g_2\in G$ and
$h,h_0,h_1\in H$, we have:

\begin{enumerate}
 \item $\d\{g_1,g_2\} = [g_1,g_2]$
 \item $\{\d h_1,\d h_2\} =[h_1,h_2]$
 \item $\{\delta h, g\}= h^{-1}(h^g)$
 \item $\{g,\delta h\}= (h^g)^{-1}h$
 \item $\{g_0,g_1g_2\}=\{g_0,g_2\}\{g_0,g_1\}^{g_2}$
 \item $\{g_0g_1,g_2\}=\{g_0,g_2\}^{g_1}\{g_1,g_2\}$
 \item $\{g_1,g_2\}\{g_2,g_1\}=1$
 \item $\{g,g\}=1$.
\end{enumerate}
In view of axiom (1), $\{-,-\}$ is called a {\it lift} of the
commutator map to $H$. By an abuse of notation, we denote a crossed
module (or a strictly stable crossed module) by $H\ra G$.
\end{defe}

\begin{defe} Let $\fG:=(H\ra G)$ and $\fG':=(H'\ra G')$ be
crossed modules. A morphism of crossed modules is a pair of
homomorphisms $a:H\ra H'$ and $b:G\ra G'$ such that the following
diagram commutes
\[
\xymatrix{
  H \ar[d]_{a} \ar[r] & G \ar[d]^{b} \\
  H' \ar[r] & G',   }
\]
and for every $h\in H$ and $g\in G$,  $a(h^g)=a(h)^{b(g)}$.

Suppose that $\fG$ and $\fG'$ are strictly stable crossed modules,
and let $\{-,-\}$ and $\{-,-\}'$ be their respective lift of the
commutator map. Then $(a,b)$ is a morphism of strictly stable
crossed modules if the following diagram commutes:
\[
\xymatrix{
  G\times G \ar[d]_{a} \ar[r]^{\quad \{-,-\}} & H \ar[d]^{b} \\
  G'\times G' \ar[r]^{\quad \{-,-\}'} & H'.   }
\]
\end{defe}

\subsection{Relationship between gr-categories and crossed modules}
\label{s:relationShipCrosModGrCat}

 \underline{Crossed modules to
gr-categories:} The passage from crossed modules to gr-categories
uses the notion of quotient groupoid.

\begin{defe}\label{d:quotientGroupoid}
Let $\Gamma$ be a group acting on a set $X$. The {\it quotient
groupoid} $[X/\Gamma]$ is the groupoid whose objects are the
elements of $X$. An arrow $x\ra x'$ in $[X/\Gamma]$ is an element
$\gamma\in \Gamma$ such that $\gamma.x=x'$.
\end{defe}

Let  $H\ra G$ be a crossed module. Then $H$ acts on $G$ by right
translation. One can show that $[G/H]$ is a gr-category.
Furthermore, a morphism of crossed modules defines a {\it strict}
1-morphism between the corresponding gr-categories. These facts
remain true if we replace crossed modules by strictly stable crossed
modules and gr-categories by strictly commutative Picard groupoids;
see, for instance, ~\cite{Cond},
\\~\cite{Br}, ~\cite{BrBook}, and ~\cite{behrang}.

\underline{Gr-categories to crossed modules:} Let $\cM$ be a
gr-category. Let $\cM'$ be a strict model of $\cM$; that is to say,
$\cM'$ is a strict gr-category which is equivalent to $\cM$. Let $G$
be the group of objects of $\cM'$ and let $H$ be the group of arrows
whose source is $\bone$. There is a natural morphism $H\ra G$ given
by taking an arrow to its target. Furthermore, $G$ acts on $H$ by
conjugation. One can show that $H\ra G$ is a crossed module and
$[G/H]$ is equivalent to $\cM$; see, for instance, ~\cite{behrang}.

Our goal is to prove a first isomorphism theorem for gr-categories.
Recall that the first isomorphism from group theory states that
every homomorphism $f:\Gam \ra M$ has a canonical factorization
$\Gam\onto \Gam/(f^{-1}(1))\inj M$.
\begin{ntn} Let $F: \cC\ra \cD$ be a functor between groupoids. Let
$d\in \cD$. Let $F^{-1}(d)$  denote the following groupoid:

 objects =
objects of $\cC$ which map to $d$

morphisms = morphisms of $\cC$ which map to $\id_d$.
\end{ntn}

Let $G$ be a group and let $\cM=(\cM,\otimes, \bone)$ be a
gr-category. Let $F:G\ra \cM$ be a monoidal functor (Example
\ref{ex:unG}). Let $H:=F^{-1}(\bone)$.
\begin{lemma}(First Isomorphism Theorem)
\label{l:firstIsoForgrStacks}
 \begin{enumerate}
 \item $H\ra G$ is a crossed module.
 \item $F$ has
a factorization $G\onto [G/H]\inj \cM$.
\end{enumerate}
\end{lemma}

\begin{proof}
By ~\cite{behrangWeakMaps}, thm. 7.10, $F$ has a factorization
$G\rar{F'} \cM'\rar{e} \cM$ where $\cM'$ is a strict gr-category,
$F'$ is a strict monoidal functor, and $e$ is an isomorphism of
gr-category. Let $H':=F'^{-1}(d)$. Then, $e$ defines an isomorphism
of gr-categories $H'\isom H$. Therefore, we may assume that $\cM$
and $F$ are strict.

By definition, the objects of $H$ are pairs $(x,\ph)$ where $x\in G$
and $\ph$ is an isomorphism $\bone \ra F(x)$ in $\cM$. Define a
multiplication on $H$ by
\[
(x,\ph).(y,\psi):=(xy, \ph\otimes \psi).
\]
Note that $\ph\otimes \psi$ is an isomorphism $\bone \ra F(x)\otimes
F(y) = F(xy)$; hence, this product is well-defined. As $\cM$ is
strictly associative, this product is associative. It is easy to
check that the above multiplication makes $H$ into a group, where
the unit is $(1,\id_{\bone_U})$. Furthermore, the map $H\ra G$,
defined by $(x,\ph)\mapsto x$, is a group morphism.

Next, define an action of $G$ on $H$ by
\[
(x,\ph)^g:=(g^{-1}xg, \id_{F(g)^{-1}} \otimes \ph \otimes
\id_{F(g)}).
\]
One checks that this map is indeed an action, making $H\ra G$ into a
crossed module. Finally, define a functor $K: [G/H]\ra \cM$ as
follows:

 $K(g):=F(g)$ for $g\in G$

$K(g\rar{h} g\delta{(h)}):=\id_g \otimes (\bone \rar{h}
\delta{(h)})$.\\
 $K$ is a monomorphism of gr-categories, giving rise to the
required factorization.
\end{proof}

\begin{rem} \label{r:strictlyStableIfPicard}
In the previous lemma, if $\cM$ were a strictly commutative Picard
groupoid, then $H\ra G$ would be a strictly stable crossed module.
\end{rem}

\subsection{Gr-stacks and Picard stacks}
\label{s:grStPicardSt}

Let $k$ be a field.
Let $k_\et$ denote the category of sheaves on the category of
schemes over $k$, equipped with the \e tale topology. The 2-category
of stacks (in groupoids) on $k_\et$ is defined in
~\cite{GiraudNonAbelianCoh}. Roughly speaking, a stack $\cS$ is a
sheaf of groupoids: for every scheme $U$ over $k$, the ``sections"
of $\cS$ above $U$ form a groupoid, which is denoted by $\cS_U$. The
gluing conditions for these groupoids is best expressed using fibred
categories and descent ~\cite{GR}, ~\cite{GiraudNonAbelianCoh}. In
this text, we will only be concerned with Deligne-Mumford stacks
~\cite{DM}.

\begin{exam} Let $G$ be an algebraic group acting on an algebraic
variety $X$. Suppose the stabilizer of every point of $X$ is an \e
tale subgroup of $G$. Then the quotient stack $[X/G]$\footnote{For
every scheme $U$ over $k$, $[X/G]_U$ equals the quotient groupoid
$[\Hom(U,X)/\Hom(U,G)]$, see Definition \ref{d:quotientGroupoid}.}
is a Deligne-Mumford stack.
\end{exam}

Gr-stacks (resp. Picard stack) are studied in ~\cite{Br}, \S 3
(resp. ~\cite{DeligneSGA4III}). Roughly speaking, a gr-stack (resp.
Picard stack) is a sheaf of gr-categories (resp. Picard groupoids).
It is easy to see that gr-stacks and Picard stacks form a
2-category. We omit the obvious analogue of the relationship between
crossed modules and gr-stacks; see, for instance, ~\cite{Br} and
~\cite{BrBook}. The only result we use in the main body of the text
is the First Isomorphism Theorem for Deligne-Mumford Picard stacks.

Let $G$ be a connected algebraic group over $k$. (Note that $G$ can
also be considered as a Deligne-Mumford gr-stack on $k$.) Let $\cP$
a Deligne-Mumford strictly commutative Picard stack and let $F:G\ra
\cM$ be a 1-morphism of gr-stacks. Let $H=F^{-1}(\bone)$.

\begin{lemma}\label{l:FirstIsoForPicardStacks}
\begin{enumerate}
\item[(i)] $H$ is an algebraic group; that is to say, it is a reduced scheme of finite
type over $k$.
 \item[(ii)]  $H\ra G$ is a strictly stable crossed module.
 \item[(iii)] $F$ has a factorization $G\onto [G/H]\inj \cP$.
 \end{enumerate}
\end{lemma}
\begin{proof}
For every scheme $U$ over $\spec(k)$, and every object $x\in H_U$,
we have $\Aut_{H_U}(x)=\{\id\}$. By ~\cite{LauMore}, cor. 8.1.1
(iii), $H$ is representable (by an algebraic space). As $H$ has a
group structure, we conclude that it is represented by a group
scheme over $k$. Moreover, the kernel of the canonical morphism
$H\ra G$ equals $\pi_1(\cP)$ and is, therefore, an \e tale group
scheme over $k$. It follows that $H$ is an algebraic group over $k$;
i.e., it is reduced of finite-type.

Statements (ii) and (iii) follow from Lemma
\ref{l:firstIsoForgrStacks} and Remark
\ref{r:strictlyStableIfPicard}.
\end{proof}

\section{Computing fundamental group schemes}
\label{a:ComputingFundamentalGroupSchemes}

\subsection{Characteristic zero}
Let $G$ be a connected algebraic over a field of characteristic
zero. The K\"{u}nneth formula implies that the two homomorphisms in
(\ref{dia:coequalizer}) are equal. It follows that
$\pi_1(G)=\pet(G)$.

\subsection{Semisimple groups}
\label{s:petSSGroup} Let $G$ be a connected semisimple group over an
algebraically closed field. Let $\pi_1^\ss(G)$ denote the weight
lattice modulo the root lattice.\footnote{The choice of the maximal
torus does not matter.} We claim that $\pet(G)=\pi_1^\ss(G)$.

 It is enough to show that a connected simply
connected semisimple algebraic group $G$ does not admit a nontrivial
central extension by finite groups. Suppose $1\ra A\ra \tG\rar{p}
G\ra 1$ is a central extension of $G$ where $A$ is a finite abelian
group, and $\tG$ is connected. It is clear that $\tG$ is semisimple.
Let $\tT$ be a connected maximal torus of $\tG$, and let $T:=p(\tT)$
denote the corresponding maximal torus in $G$. Since the center of
$\tG$ is a subgroup of $\tT$, we have a central extension $1\ra A\ra
\tT\ra T\ra 1$. This, in turn, defines an inclusion
\[
X(T)\inj X(\tT) \inj \Lambda(T)=\Lambda(\tT),
\]
where $X$ (resp. $\Lambda$) denotes the root (resp. weight) lattice.
As $G$ is simply connected, $X(T)=\Lambda(T)$; thus, $X(T)=X(\tT)$,
implying that $A$ is trivial.

\subsection{Finite fields} Let $\Fq$ denote a finite field with
$q$ elements. Let $G$ be a connected algebraic group over $\Fq$.
 Let $\Fr$ denote the Frobenius automorphism $x\mapsto
x^q$. Let $A$ be a finite abelian group (considered as a discrete
group scheme over $\Fq$). A central extension
\[
1\ra A\ra \tG\rar{\pi} G\ra 1
\]
 defines a homomorphism\footnote{This is an instance of Grothendieck's sheaf-function
 correspondence. The sheaf in question, is the sheaf of local sections of $\pi$.}
  $f_\pi: G(\Fq)\ra A$ as follows. Given $g\in G(\Fq)$,
pick $g'\in \pi^{-1}(g)$. Note that $\Fr(g')=ag'$ for some $a\in A$.

\begin{lemma} (~\cite{MityaThesis}, \S 1.4.) \label{l:pdiscCommutative}
The map $\H^2(G,A)\ra \Hom(G(\Fq),A)$ given by the above
construction is an injection. If $G$ is commutative, then it is an
isomorphism.
\end{lemma}

\begin{cor} Let $G$ be a connected commutative algebraic group over
$\Fq$. Then $\pdisc(G)=G(\Fq)$.
\end{cor}

\subsection{Additive group}
Let $k$ be an algebraically closed field of positive characteristic.
Let $G$ denote the additive group over $k$. Our aim is to prove
\[
\pet(G)=\Hom(k,\qz)
\]
\begin{rem}[Pontryagin Duality] The category of profinite abelian
groups is anti-equivalent to that of discrete torsion abelian
groups: to a profinite abelian group $\Pi$, one associates
$\Gamma:=\Hom_\cont (\Pi,\qz)$. The inverse functor is $\Gamma \ra
\Hom(\Gamma,\qz)$.
\end{rem}

By Pontryagin duality, it is enough to show that $\H^2(G,\qz)=k$. As
$G$ is connected and commutative, the natural injection
$\Ext(G,\qz)\inj \H^2(G,\qz)$
 is an isomorphism. Thus, it is enough
to show that \newline
$\Ext(G,\qz)=k$. On the other hand, $G$ is
killed by multiplication by $p$; hence, the long exact sequence
corresponding to
\[
0\ra p\bZ/\bZ\ra \qz\rar{p}\qz\ra 0
\]
shows that $\Ext(G,\qz)=\Ext(G,\Fp)$. Thus, we are reduced to
proving $\Ext(G,\Fp)=k$.

Observe that
\[
\Ext(G,\Fp)\subseteq
\H_\et^1(G,\Fp)=\H_\et^1(\bA_k^1,\Fp)=k[x]/{A(k[x])},
\]
 where $A(u)=u^p-u$. Furthermore, every element of $\H_\et^1(G,\Fp)$ can be
written as
\[
h=\displaystyle \sum_{p\nmid i} c_i\overline{x^i},
\]
 where
$\overline{x^i}$ is the image of $x^i$ in the quotient
$k[x]/{A(k[x])}$. By Lemma \ref{l:petQuitient}, $h\in \Ext(G,\Fp)$
if and only if the polynomials
\[
\sum_{p\nmid i} c_i(x^i+y^i), \quad \quad \textrm{and} \quad \quad
\sum_{p \nmid i} c_i(x+y)^i\notag
\]
represent the same element of $\H_\et^1(\bA^2,
\Fp)=k[x,y]/{A(k[x,y])}$. This happens if and only if the polynomial
$t(x,y)=\displaystyle \sum_{p\nmid i}c_i[ (x+y)^i-x^i-y^i]$ is of
the form $u^p-u$. This is possible if and only if $c_i=0$ for all
$i\neq 1$ (otherwise, the degree of $t$ is not divisible by $p$).
Therefore, $\Ext(G,\Fp)$ identifies with
\[
\{c\overline{x}\, |\, c\in k \}\subseteq
k[x]/{A(k[x])}=\H_\et^1(\bA^1,\Fp).
\]

\bibliography{fundamentalGroupYo}
\end{document}